\documentclass[11pt]{article}
\usepackage{amsfonts}
\begin{document}

\begin{center}
\textbf{About one class polynomial problems with not polynomial certificates}

Kochkarev B.S.

{\it University of Kazan, Russia}

E-mail: bkochkarev@rambler.ru
\end{center}

Abstract. We build a class of polynomial problems with not polynomial certificates. The parameter concerning which are defined efficiency of corresponding algorithms is the number $n$ of elements of the set has used at construction of combinatory objects (families of subsets) with necessary properties.
\vskip3mm
Let $\Sigma =\{0,1\}$ is set of two elements. An language $L$ over $\Sigma$ is any set of strings made up of symbols from $\Sigma$. We denote the empty string by $e$, and the empty language by $\emptyset$. The language of all strings over $\Sigma$ is denoted $\Sigma^*$. Every language $L$ over $\Sigma$ is a subset of $\Sigma^*$. Let $n$ is a natural number and $S$ is the ordered set $(a_{1},a_{2},\ldots,a_{n})$. In the capacity of combinatorial objects we will to consider the families of subsets of the set $S$ with necessary properties. Any subset $A$ of $S$ can be present in look of string from $n$ elements $(\sigma_{1},\sigma_{2},\ldots,\sigma_{n})$, where $\sigma_{i} =1$, if $a_{i}\in A$ and $\sigma_{i} =0$, if $a_{i}\notin A$. In the capacity of languages we will consider the families of subsets of the set $S$ with definite properties. Thus a language in our case is a subset of the strings $(\sigma_{1},\sigma_{2},\ldots,\sigma_{n})$.

In 1964  Alan Kobham [1] and, independently, in 1965  Jack Edmonds [2] have entered a concept of complexity class $P$.

Definition 1 [1,2]. A language $L$ belong to $P$ if there is an algorithm $A$  that decide $L$ in polynomial time $(\leq O(n^{k}))$ for a constant $k$.

Class of problems $P$ is called polynomial.

According to [3] J. Edmonds has entered also the complexity class $NP$. This is the class of problems (languages) that can be verified by a polynomial-time algorithm.

 Definition 2 [3]. A language $L$ belongs to $NP$ if there exists a two-input polynomial-time algorithm $A$ and such polynomial $p(x)$ with whole coefficients that

$L=\{x\in \{0,1\}^{n}:  there\,  exists\,  a\,  certificate\,  $y$\, with \mid y \mid\leq p(\mid x\mid)\, and\, A(x,y)=1\}$

 In this case we say that the algorithm $A$ verifies language $L$ in polynomial time.

 According to definition 2 if $L\in P$ and $\mid y \mid \leq p(\mid x \mid)$ then $L\in NP$. But if $L\in P$ and length of the certificate not polynomial from length $x$ then $L\notin NP$.

 J.Edmonds has conjectured also that $P\not= NP$ which so far  is not proved.

 In 1971 S.A.Cook has put the question: "whether can the verification of correctness of the decision of a problem be more long than the decision itself independently of algorithm of verification?" This problem have a relation to cryptography. In other formulation this problem look so: "whether can be build a cipher such that his decipher algorithmically more complicated than find the cipher?"

 The RSA public-key cryptosystem is based on the dramatic difference between the ease of finding large prime numbers and the difficulty of factoring the product of two large prime numbers. In connection with it arise the question in Cook's problem.

 In 2008 [4] we have proposed a model of decision Cook's problem: let $M$ and $M'$ are two sets such that $M$ is decidable set in polynomial time, let then there exists the injective map $\phi$ of $M$ in $M'$ such that for any $m\in M$ $\phi (m)$ find in not polynomial time. In [5,6] we have cited some realizations of this model. Here we cite his yet one realization.

 Definition 3 [7]. A family $F$ of subsets of the set $S$ is called Sperner, if no element $A\in F$ is a subset of another element $A^\prime\in F$.

Definition 4 [8-10]. A Sperner family (S.f.) $F$ is called maximal, if for any $A\subset S, A\notin  F$, one can find $A^\prime \in F$ such that $A\subset A^\prime$ or $A^\prime\subset A$.

Definition 5 [ibid] We say that a S.f. $F$ has the type $(k,k+1)$, if $\mid A\mid \in \{k,k+1\}$ for any $A\in F$.

Let $F$ be a maximal Sperner family (m.S.f.) of the type $(k,k+1)$, $k\not= 0,n-1$. Thus, we do not consider the m.S.f. $\emptyset$ and $S$. Let $p_{i}$ stand for the number of elements $A\in F,\mid A\mid= k$,
 which do not contain the element $a_{i}\in S$; let $q_{i}$ stand for the number of elements   $A\in F,\mid A\mid = k+1$ which contain the element $a_{i}$. Let $r_{i} =p_{i} + q_{i}$, $r = max r_{i}, i = \overline{1,n}.$Evidently, with any $n \geq 3$ the following inequality is true $r_{i} \leq { n-1 \choose k }.$
It is well-known [8] that studying m.S.f. of the type $(k,k+1)$, it suffices to consider the case $k \leq \lfloor \frac {n}{2}\rfloor.$

If $F$ is a S.f. of the type $(k,k+1)$, then $F^{(k)}, F^{(k+1)}$ stand, correspondingly, for families of subsets $A\in F, \mid A\mid = k, A^\prime \in F, \mid A^\prime\mid = k+1.$

If $F^{(k +1)}$ is family such that for any $A\in F^{(k +1)}$ $a_{i}\in A$ then we denote by $F^{(k +1)}\backslash \{a_{i}\}$ the S.f. of the set $S\backslash \{a_{i}\}$ obtained by expel of each subset $A\in F^{(k +1)}$ of element $a_{i}$.

Theorem 1 [9]. If $F$ m.S.f. then

              ${n-1 \choose k } \leq \mid F\mid \leq {n \choose k +1}$

Definition 6. S.f. $F^{(k +1)}$ is called the admissible fragment of m.S.f. $F^\prime \supset F^{(k +1)}$ if there is $i\in \overline {1,n}$ such that for any $A\in F^{(k +1)}$ $a_{i}\in A$.

Theorem 2 [6]. There exists the injective map of set admissible fragments $\{F^{(k +1)}\}$ in set of m.S.f. $\{F^\prime\}$ with $r(F^\prime) = {n-1 \choose k }.$

Proof. Let $F^{(k +1)}$ is an admissible fragment such that for any $A\in F^{(k +1)}$ $a_{i}\in A.$ Then we denote by $F$ S.f. $F^{(k +1)}\setminus \{a_{i}\}.$ Let then $F_{1}$ is $\{A, \mid A\mid = k, a_{i}\notin A\}\setminus F.$ Evidently family

$F^\prime = F^{(k +1)}\bigcup F_{1}\bigcup G\bigcup G^\prime,$

where $G^\prime$ is family of subsets $A,\mid A\mid = k, a_{i}\in A$ , incomparable by inclusion with elements of $F^{(k +1)}$ and $G$ is family of subset $A^\prime, \mid A^\prime \mid = k+1, a_{i}\notin A^\prime,$
incomparable by inclusion with elements of $F_{1}$ is the family corresponding the admissible fragment $F^{(k +1)}.$ Evidently $r(F^\prime) ={n-1 \choose k}$. The map receiving is injective by construction.

Let $n$ is enough large odd number. We consider S.f. $F^{(\lceil \frac {n}{2}\rceil)}$ of subsets of the set $S$ such that for any $A\in F^{(\lceil \frac {n}{2}\rceil)}$, $a_{i}\in A$, if $i\in\{1,2, \ldots, \lceil \frac {n}{2}\rceil - l \}$ where $l$ is a constant. Evidently $\mid F^{(\lceil \frac {n}{2}\rceil)}\mid$ is polynomial of $n$.

Theorem 3. Family $F^{(\lceil \frac {n}{2}\rceil )}$ is decided in polynomial-time.

Proof. Really since $\mid F^{(\lceil \frac {n}{2}\rceil)}\mid$ is the polynomial from $n$ there exists polynomial algorithm $T$ which in polynomial number of steps definit for any subset $A \subset S$ belongs it $F^{(\lceil \frac {n}{2}\rceil)}$ or no.

We consider two problems: straight (to find $A\in F^{(\lceil \frac {n}{2}\rceil )}$) and reverse (by finding $A$ to build the m.S.f. $F^\prime$ corresponding $A$). In this connection we suppose $\{A\}$ as the set of ciphers and $\{F^\prime\}$ as the corresponding deciphers. Since $\mid F^\prime \mid > {n-1 \choose  \frac {n-1}{2}}$  is exponent [11] then decipher algorithmically more complex than find the cipher. Thus, at check of correctness of decision the straight problem the length of certificate is not polynomial. From here our  straight problem of $P$ does not belong to $NP$, that is $P\not= NP$.We note at last in addition to [12] once more that the statement $P\subseteq NP$ is error.

\centerline{REFERENCES}
\smallskip

1. Cobham A. The intrinsic computational difficulty of functions // In Procedings of the 1964 Congress for Logic, Methodology, and the Philosophy of Science.- North-Holland, 1964.-P.24-30.

2. Edmonds J. Paths, trees and flowers // Canadian Journal of Mathematics.-1965-Vol.17.-P.449-467.

3. Cormen T.H., Leiserson Ch.E., Rivest R.L., Introduction to Algorithms, MIT Press, 1990. 

4. Kochkarev B.S., On Cook's problem, www.math.nsc.ru//simalglog/ses 2008e.html

5. Kochkarev B.S., Prilojenie monotonnykh funktsij algebry logiki k probleme Kuka, Nauka v Vuzakh: matematika, fizika, informatika, Tezisy dokladov Mejdunarodnoj nauchno-obrazovatelnoj konferentsii,2009,pp.274-275.

6. Kochkarev B.S.,K probleme Kuka, Matematicheskoje obrasovanije v shkole i v vuze v uslovijakh perekhoda na novye obrazovatelnye standarty, Materialy Vserossijskoj nauchno-practicheskoj konferentsii s mejdunarodnym uchastiem, Kazan, pp.133-136, 2010.

7. Sperner E., Ein Satz  \" uber Untermengen einer edlichen Menge, Math.Z. 27 (1928),pp.544-548.

8. Kochkarev B.S., Structural Properties of One Class of Maximal Sperner Families of Subsets of a Finite Set, In Proceedings of International Conference, "Logic and Applications" on the occasion of the 60th anniversary of the Academician Yu.l. Erchov, Novosibirsk, pp.61-62, 2000.

9. Kochkarev B.S., Structure Properties of a Certain Class of Maximal Sperner Families of Subsets, Russian Mathematics (Iz.VUZ) 59 (7), pp.35-40, 2005.

10. Kochkarev B.S., Admissible Values of One Parameter for Maximal Sperner Families of Subsets of the Type (k,k+1), Russian Mathematics (Iz. VUZ), Vol 52, 6 pp.22-24,2008.

11. Yablonskij S.V. Vvedenie v diskretnuju matematiku, 384, 1986.

12. Kochkarev B.S. Gipoteza J.Edmondsa i problema S.A.Kuka, Vestnik TGGPU, (24),№2, pp.23-24, 2011.
\end{document}